\documentclass[a4paper,reqno,twoside]{amsart}
\usepackage{amsmath,amsfonts,amssymb,enumerate,url}
\usepackage[english]{babel}
\usepackage[all]{xy}

\theoremstyle{plain}
 \newtheorem{theorem}{Theorem}[section]
 \newtheorem{tthr}[theorem]{Theorem}
 \newtheorem{corr}[theorem]{Corollary}
 \newtheorem{lmr}[theorem]{Lemma}
 \newtheorem{prr}[theorem]{Proposition}
 
\theoremstyle{definition}
 \newtheorem{defnr}[theorem]{Definition}
 \newtheorem{remnr}[theorem]{Remark}

\def\brm{\begin{remnr}}
\def\erm{\end{remnr}}
\def\bdr{\begin{defnr}}
\def\edr{\end{defnr}}
\def\bp{\begin{proof}}
\def\ep{\end{proof}}
\def\btr{\begin{tthr}}
\def\etr{\end{tthr}}
\def\bpr{\begin{prr}}
\def\epr{\end{prr}}
\def\bcr{\begin{corr}}
\def\ecr{\end{corr}}
\def\blr{\begin{lmr}}
\def\elr{\end{lmr}}
\def\beq{\begin{equation}}
\def\eeq{\end{equation}}
\def\bcs{\begin{cases}}
\def\ecs{\end{cases}}

\def\dt{,\dots,}
\def\dpt{+\dots+}
\def\ind{{\rm ind}}
\def\Bd#1{\mathbb #1}
\def\quat#1{\textquotedblleft #1\textquotedblright}

\def\a{\ar@{-}}
\def\*{\bullet}

\title{TWO COMBINATORIAL FORMULAS CONCERNING MARKED PARTITIONS}
\author{F.~V.~Weinstein}
\address
{Giacomettistrasse 33A, CH-3006 Bern, Switzerland.}
\email{felix.weinstein46@gmail.com}

\begin{document}

\begin{abstract}
A partition of degree $n$ is a decomposition
$n=i_1+i_2+\dots+i_q$, where ${i_1,i_2,\dots,i_q}$ are positive integers
called the parts of the partition.
Let $\lambda>0$ be an integer.
The partition is said to be a $\lambda$--partition
if $i_{a+1}-i_a\geqslant \lambda$ for all $a$ such that $1\leqslant a<q$.

The main result of this note are combinatorial formulas,
which express the quantity of $1$-partitions of a given degree
in terms of the $\lambda$--partitions of the same degree, where $\lambda=2$ or $\lambda=3$,
some special parts of which are marked depending on $\lambda$.
The presented proofs of both formulas are bijective.

It is shown that for $\lambda=3$ the corresponding formula
is equivalent to the classical Sylvester identity.

The obtained combinatorial formulas as well as their bijective proofs are generalized
to the quantities of $1$--partitions, all parts of which are $\geqslant k$ for any fixed integer $k\geqslant 1$.
\end{abstract}

\maketitle
\markboth{TWO COMBINATORIAL FORMULAS CONCERNING MARKED PARTITIONS}{F.~V.~WEINSTEIN}

\section*{\bf Introduction}

The main result of this note are two new combinatorial formulas.
Each of them express the quantity of $1$--partitions of degree $n$ in terms of $\lambda$--partitions,
where either $\lambda=2$, or $\lambda=3$, some parts of which are marked
depending on $\lambda$. Our proofs of these formulas are bijective.

Using generating functions both formulas simultaneously can be written as an
identity between formal power series.
It turns out that for $\lambda=3$ this identity is equivalent to
the classical identity
\beq\label{Syl}
\prod^\infty_{r=1}(1+tx^r)=
1+\sum_{q=1}^\infty
\,t^qx^{\frac{3q^2-q}{2}}\;
\frac{(1+tx)\dots(1+tx^{q-1})(1+tx^{2q})}
{(1-x)\dots(1-x^{q-1})(1-x^q)}\;,
\eeq
obtained by Sylvester in 1882 in article \cite{MR1505328} (p.282).
The known proofs of it usually use a machinery of
generating functions combined with additional combinatorial arguments.
Such is the original proof by Sylvester and the proof from book \cite{MR1634067} (Th.9.2).
Our approach provides the Sylvester identity with a clear combinatorial interpretation,
as well as with a bijective proof of it.

For $\lambda=2$ the corresponding identity
does not have such a nice form as for $\lambda=3$.
However, it implies the identity
\beq\label{Syl0}
\prod_{r=1}^\infty(1-x^{2r-1})=
1+\sum_{q=1}^\infty\frac{(-1)^qx^{q^2}}{(1-x^2)(1-x^4)\dots(1-x^{2q})}\;.
\eeq

The note is organized as follows.
In \S\ref{sec1}, the necessary notations and definitions are introduced.
In \S\ref{sec002}, the main result is formulated (Theorem \ref{t_main})
and some of its  corollaries are presented.
In particular, it is shown how it implies formulas \eqref{Syl} and \eqref{Syl0}.
In \S\ref{sec2} and \S\ref{sec03}, for $\lambda=3$ and $\lambda=2$,
respectively, the mentioned bijective correspondences are constructed.
These bijections
are quite different for $\lambda=3$ and for $\lambda=2$ (however, see Remark \ref{rem2}).
In \S\ref{sec04} the obtained combinatorial formulas
as well as their bijective proofs are generalized to
the sets of $1$--partitions with all parts $\geqslant k$ for any fixed integer $k\geqslant 1$ (Theorem \ref{t_maink}).

\section{\bf Marked $\lambda$--partitions}\label{sec1}

In what follows we use the notation:
\smallskip

\qquad $|M|$\; is the cardinality of the finite set $M$.
\smallskip

\qquad $[a,b]:=\{r\mid a\leqslant r\leqslant b\}$, where $a,b,r$ are integers.

\bdr
A \emph{partition} is a finite set of positive integers $I=\langle i_1,i_2,\dots,i_q\rangle$,
called its \emph{parts}.
The numbers
\[
\|I\|:=i_1\dpt i_q\qquad\text{and}\qquad |I|=q
\]
are called the \emph{degree} and the \emph{length} of partition $I$, respectively.
\edr

\bdr
A \emph{marked partition} is a pair $\langle I;J\rangle$,
where $I$ is a partition and $J\subset I$.
The parts of $I$ belonging to $J$ are called the \emph{marked parts} of $\langle I;J\rangle$.
The numbers
\[
\left\|\langle I;J\rangle\right\|:=\|I\|\qquad\text{and}\qquad\left|\langle I;J\rangle\right|:=|I|+|J|
\]
are called the \emph{degree} and the \emph{length} of the marked partition $\langle I;J\rangle$, respectively.

Any partition $I$ is interpreted as the marked partition $\langle I;\emptyset\rangle$.
Define $\langle I_1;J_1\rangle\cup\langle I_1;J_1\rangle:=\langle I_1\cup I_2;J_1\cup J_2\rangle$.
\edr

Instead of separately indicating the set of marked parts, we often underline them:
$\big\langle 1,5,8\,;\,5\big\rangle=\langle 1,\underline{5},8\rangle$.

\bdr
Let $\lambda>0$ be an integer.
A \emph{$\lambda$-partition} is a pair $(\lambda,I)$,
where $I=\langle i_1,i_2\dt i_q\rangle$ is a partition
such that $i_{a+1}-i_a\geqslant \lambda$ for any $a\in[1,q-1]$.

We say that \emph{marked partition $\langle I;J\rangle$ is a $\lambda$-partition}, if $I$ is a $\lambda$-partition.
\edr

We write $\lambda$-partitions as usual partitions, emphasising
that we only consider $\lambda$-partitions.
For example, one may treat $\langle 2,5,8\rangle$ as a $1$-, $2$-, or $3$-partition.
These objects are not the same.

\bdr
A \emph{dense $\lambda$-partition} is a partition $\langle i_1,i_2\dt i_q\rangle$,
where $i_{a+1}-i_a=\lambda$ for any $a\in[1,q-1]$.
\edr

For $\lambda$--partitions $I_1,I_2$ such that $\min(I_2)-\max(I_1)>\lambda$, we write
the $\lambda$--partition $I_1\cup I_2$ as $I_1\sqcup I_2$.
\bdr
A \emph{canonical form of the $\lambda$-partition $I$} is a decomposition
$I=I_1\sqcup I_2\sqcup\dots\sqcup I_s$, where
$\lambda$-partitions $I_1,I_2,\dots,I_s$ are dense.
\edr

{\it In what follows we assume that $\lambda=2$ or $\lambda=3$.}

\bdr
Let
$I=I_1\sqcup I_2\sqcup\dots\sqcup I_s$ be the canonical form of the $\lambda$--partition $I$.
Define
\[
\ind_\lambda(I):=\ind_\lambda(I_1)+\ind_\lambda(I_2)+\dots+\ind_\lambda(I_s),
\]
where for a dense $\lambda$--partition $I=\langle i,i+\lambda\dt i+\lambda(q-1)\rangle$ we define
\[
\ind_2(I):=
\begin{cases}
0&\text{if $i=1$ or $i\equiv 0\mod 2$},\\
1 &\text{otherwise},
\end{cases}
\qquad\qquad
\ind_3(I):=
\begin{cases}
0&\text{if $i=1,2$},\\
1&\text{if $i>2$}.
\end{cases}
\]
The number $\ind_\lambda(I)$ is called the \emph{index} of the $\lambda$--partition $I$.

The minimal part of $I_a$ for any $a\in[1,s]$,
is called a \emph{leading part} of $I$ if $\ind_\lambda(I_a)=1$.
Thus, $\ind_\lambda(I)$ is the number of the leading parts of the $\lambda$--partition $I$.
\edr

\bdr
A marked $\lambda$--partition $\langle I;J\rangle$ is called \emph{regular},
if $J$ is a subset of the leading parts of $I$.
\edr
\smallskip

{\it Examples}: The canonical form of the $2$-partition $I=\langle 1,3,9,11,14\rangle$
is $I=\langle 1,3\rangle\sqcup\langle 9,11\rangle\sqcup\langle 14\rangle$.
Since $9$ is its single leading part, $\ind_2(I)=1$.

The canonical form of the $3$-partition $I=\langle 2,5,8,12,15,19\rangle$
is $I=\langle 2,5,8\rangle\sqcup\langle 12,15\rangle\sqcup\langle 19\rangle$.
The set of leading parts of $I$ consists of parts $12$ and $19$. Therefore, $\ind_3(I)=2$.

\section{\bf Main result and its corollaries}\label{sec002}

In what follows we use the notation:
\[
\begin{array}{ll}
D(n,q)&\text{is the set of 1--partitions of degree $n$ and of length $q$}.\\[.5mm]
D(n)&\text{is the set of 1--partitions of degree $n$}.\\[.5mm]
M_\lambda(n,q)&\text{is the set of regular marked $\lambda$--partitions of degree $n$ and of length $q$}.
\end{array}
\]

The main result of this note is the following

\btr\label{t_main}
For $\lambda=2$ or $\lambda=3$, we have $|D(n,q)|=|M_{\lambda}(n,q)|$. In particular,
\beq\label{main1}
\prod_{r=1}^\infty(1+tx^r)=1+\sum_{q=1}^\infty\;
\sum_{n=\frac{q(q+1)}{2}}^\infty\left|M_{\lambda}(n,q)\right|\,t^q\,x^n.
\eeq
\etr

For example, each of the following sets contains 7 elements:
\[
D(12,3)=
\left\{
\langle 1,2,9\rangle,\langle 1,3,8\rangle,
\langle 1,4,7\rangle,\langle 1,5,6\rangle,
\langle 2,3,7\rangle,\langle 2,4,6\rangle,\langle 3,4,5\rangle
\right\},
\]
\[
M_2(12,3)=
\left\{
\langle 1,\underline{11}\rangle,\langle\underline{3},9\rangle,\langle 3,\underline{9}\rangle,
\langle\underline{5},7\rangle,\langle 1,3,8\rangle,\langle 1,4,7\rangle,\langle 2,4,6\rangle\right\},
\]
\[
M_3(12,3)=
\left\{
\langle 1,\underline{11}\rangle,\langle 2,\underline{10}\rangle,\langle\underline{3},9\rangle,\langle 3,\underline{9}\rangle,
\langle\underline{4},8\rangle,\langle 4,\underline{8}\rangle,\langle 1,4,7\rangle
\right\}.
\]

To obtain some corollaries of Theorem \ref{t_main} it is convenient to present it in a more detailed form.

Namely, let $I$ be a $\lambda$--partition with $\ind_\lambda(I)=\alpha$.
The quantity of regular $\lambda$--partitions ${\langle I;J\rangle}$
with $|J|=m$ is equal to $\binom{\alpha}{m}$.
Thus, the binomial formula and identity $|D(n,q)|=|M_{\lambda}(n,q)|$ implies that
\beq\label{main}
\sum_{q=1}^\infty|D(n,q)|t^q=\sum_{q=1}^\infty\sum_{\alpha=0}^q|R_\lambda(n,q;\alpha)|\,(1+t)^\alpha t^q,
\eeq
where $R_\lambda(n,q;\alpha)$ denotes the set of
$\lambda$--partitions of degree $n$, of length $q$, and of index $\alpha$
(by definition $(1+t)^0=1$ for any $t$ including $t=-1$).
Note that the sums in each side of equality \eqref{main} are finite.

For $t=1$, the equality \eqref{main} turns into the expression
\[
|D(n)|=\sum_{\alpha=0}^\infty|R_\lambda(n;\alpha)|\;2^\alpha.
\]
where $R_\lambda(n;\alpha)$ denotes the quantity of $\lambda$--partitions
of degree $n$ and of index $\alpha$.
For example,
\begin{gather*}
|D(12)|=\sum_{\alpha=0}^\infty|R_2(12;\alpha)|\;2^\alpha=5\cdot 2^0+3\cdot 2^1+1\cdot 2^2=15,\\
|D(12)|=\sum_{\alpha=0}^\infty|R_3(12;\alpha)|\;2^\alpha=1\cdot 2^0+3\cdot 2^1+2\cdot 2^2=15.
\end{gather*}
Since any $3$-partition of index $0$ and length $q\geqslant 1$ is either
$\langle 1,4,\dots,3q-2\rangle$, or $\langle 2,5,\dots,3q-1\rangle$, then
\[
|R_3(n,q;0)|=
\bcs
1&\text{if $n=(3q^2\pm q)/2$},\\
0&\text{otherwise}.
\ecs
\]
Therefore, for $\lambda=3$ and $t=-1$, formula \eqref{main}
implies the \emph{Euler's Pentagonal Theorem}:
\[
\prod_{k=1}^\infty(1-x^k)=
1+\sum_{q=1}^\infty(-1)^q\left(x^\frac{3q^2-q}{2}+x^\frac{3q^2+q}{2}\right).
\]

Summing in both sides of formula \eqref{main} over $n$ we can present formula \eqref{main1} in the form
\beq\label{eq001}
\prod_{r=1}^\infty(1+tx^r)=1+\;\sum^\infty_{q=1}A_\lambda(x,t;q)\;t^q,
\qquad\text{where}\qquad
A_\lambda(x,t;q)=\sum^\infty_{n=1}\sum^q_{\alpha=0}|R_\lambda(n,q;\alpha)|\;(1+t)^\alpha x^n.
\eeq

For $\lambda=2$, partitions of index $0$ exist for any degree $n\neq 3$.
If $\ind_2(I)=0$, then either all parts of $I$ are even, or
$I=\langle 1,3,\dots,2q-1, 2i_1,2i_2,\dots,2i_m\rangle$,
where $q\geqslant 1$ and $i_1\geqslant q+1$.
Since ${1+3+\dots+(2q-1)=q^2}$, we obtain
\[
1+\sum_{n=1}^\infty\sum_{q=1}^\infty\;|R_2(n,q;0)|\,t^qx^n=
\left(1+\sum_{q=1}^\infty\frac{t^qx^{q^2}}{(1+tx^2)(1+tx^4)\dots(1+tx^{2q})}\right)
\prod_{s=1}^\infty(1+t\,x^{2s}).
\]
For $t=-1$, this formula together with formula \eqref{eq001} imply formula \eqref{Syl0}.

To conclude this section,
let us show that for $\lambda=3$, identity \eqref{eq001} is equivalent to the Sylvester identity.
Indeed, formula \eqref{main} implies that the coefficient
of expansion of $A_3(x,t;q)$ in power series in $t^hx^n$
is equal to the quantity of regular $3$--partitions $\langle I;J\rangle$
such that $\|I\|=n,|I|=q$, and $|J|=h\leqslant q$.

The property of partition $\langle i_1\dt i_q\rangle$ to be a 3-partition
is equivalent to the following property of its conjugate partition:
it has $i_1\geqslant 1$ parts equal to $q$,\;$i_2-i_1\geqslant 3$ parts equal to
$q-1$, and so on, $i_q-i_{q-1}\geqslant 3$ parts equal to $1$.
Therefore,
\begin{multline*}\label{A31}
A_3(x,t;q)=\left(x^q+x^{2q}+(1+t)\sum^\infty_{r=3}x^{rq}\right)
\left(x^{3(q-1)}+(1+t)\sum^\infty_{r=4} x^{r(q-1)}\right)\dots
\left(x^3+(1+t)\sum^\infty_{r=4}x^r\right)\\
=x^{\frac{3q^2-q}{2}}\;\frac {(1+tx)(1+tx^2)\dots
(1+tx^{q-1})(1+tx^{2q})} {(1-x)(1-x^2)\dots (1-x^{q-1})(1-x^q)}\,.
\end{multline*}
Indeed, if $i_1\geqslant 3$ or $i_a-i_{a-1}>3$, where $a\in[2,q]$,
the part $i_1$ or $i_a$
can be either marked (coefficient $t$), or not (coefficient $1$).
Substituting this expression into formula \eqref{eq001} gives identity \eqref{Syl}.

\begin{remnr}
Calculations in article \cite{{MR1254731}} imply an interesting formula,
which reminds of formula \eqref{main} for $\lambda=3$ and, in fact, is related to it.
Namely, for partition $I=\langle i_1\dt i_q\rangle$, define
\begin{gather*}
U(I):=\sum_{a=1}^q\binom{i_a}{3}+2\sum_{1\leqslant a<b\leqslant q}i_a i_b-3\sum_{a=1}^q(q-a)i^2_a,\qquad
V(I):=\sum_{a=1}^q\binom{i_a}{3}-\sum_{1\leqslant a<b\leqslant q}i_a i_b.
\end{gather*}
Then, for any integer $n>0$, we have
\beq\label{Lie}
\sum_{I\in D(n)}U(I)\;t^{|I|}=\sum_{N\in R_3(n)}V(N)\;(1+t)^{\ind_3(N)}t^{|N|},
\eeq
where $R_3(n)$ denotes the set of $3$--partitions of degree $n$.
The proof of identity \eqref{Lie} in \cite{{MR1254731}} uses Lie algebras.
It would be interesting to obtain a direct proof of this formula.
\end{remnr}

\section{\bf Construction of a bijection $S:D(n,q)\rightarrow M_3(n,q)$}\label{sec2}

The \emph{Ferrers diagram of the partition}
$I=\langle i_1,i_2\dt i_q\rangle$ is the set of points
$(a,b)\in\Bd{Z}\times\Bd{Z}$ (\quat{vertices}) such that $b\in[1,q]$ and $a\in[1,i_b]$.
The \emph{diagonal of the partition} $I$ is the set of vertices with $a+b=q+1$.

Let us enumerate the diagonal vertices from the bottom to the top.
Denote by $x_i$ the number of diagram vertices in the row to the right of the $i$th diagonal vertex,
including this vertex, and let $y_i$ be the number of diagram vertices located in the column
strictly below the $i$th vertex.
If the diagonal of $I$ contains $r$ vertices, then
we can interpret $I$ as a pair of integer sequences
\[
I=\left\langle x_1,x_2\dt x_r\mid y_1,y_2\dt y_r\right\rangle,\qquad\text{where\; $1\leqslant x_1<x_2<\dots<x_r,\;\;0\leqslant y_1<y_2<\dots<y_r$}.
\]
Such notation is called the \emph{Frobenius form of partition $I$}.
The next claim is obvious:
\blr\label{lms}
The partition $I=\left\langle x_1,x_2\dt x_r\mid y_1,y_2\dt y_r\right\rangle$
is a $1$--partition, if and only if
\begin{enumerate}
\item[\rm (1)] $x_{i+1}-x_i\geqslant 2$ for $i=1,2\dt r-1$.
\item[\rm (2)] $y_{i+1}-y_i=1$ or $2$ for $i=1,2\dt r-1$.
\item[\rm (3)] $y_1=0$ or $1$.
\item[\rm (4)] If $x_1=1$, then $y_1=0$.
\end{enumerate}
\elr

Let $I=\left\langle x_1,x_2\dt x_r\mid y_1,y_2\dt y_r\right\rangle$ and
let $a_1,a_2,\dots,a_s$, where
$1\leqslant a_1<a_2<\dots<a_s\leqslant r$, be all numbers such that
$y_{a_k}-y_{a_k-1}=2$, where  $y_0:=-1$.
Define
\[
R(I)=\big\langle\widehat{I};\widehat{J}\big\rangle,\quad\text{where\quad
$\widehat{I}=\left\langle x_1+y_1,x_2+y_2\dt x_r+y_r\right\rangle$,\;\;
$\widehat{J}=\left\langle x_{a_1}+y_{a_1},x_{a_2}+y_{a_2}\dt x_{a_s}+y_{a_s}\right\rangle$}.
\]
Lemma \ref{lms} implies that $R(I)\in M_3(n,q)$.
\medskip

\emph{Example}: For the partition
$I=\langle 2,3,5,6,8,9\rangle=\langle 2,4,7,9\mid 0,2,4,5\rangle$ with diagram
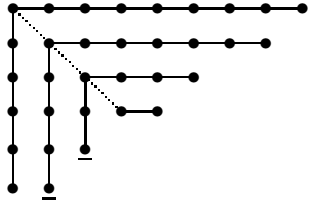
\[
\xy
\xymatrix
@=4mm @M=-0.5mm @H=-0.4mm{
\a[r]\a[d]\ar@{.}[dr]\*&\a[r]\*&\a[r]\*&\a[r]\*&\a[r]\*&\a[r]\*&\a[r]\*&\a[r]\*&\*\\
\a[d]\*&\a[d]\a[r]\ar@{.}[dr]\*&\a[r]\*&\a[r]\*&\a[r]\*&\a[r]\*&\a[r]\*&\*&\\
\a[d]\*&\a[d]\*&\a[r]\a[d]\ar@{.}[dr]\*&\a[r]\*&\a[r]\*&\*\\
\a[d]\*&\a[d]\*&\a[d]\*&\a[r]\*&\*\\
\a[d]\*&\a[d]\*&{\underline{\bullet}}\\
\*&{\underline{\bullet}}
}
\endxy
\]
the parts of $S(I)$ are the numbers of vertices connected by the solid lines.
Such a part is \emph{marked} if the lowest vertex on such a line
is located strictly below the diagonal and
there are no more vertices in the row to the right of this vertex.
Thus,
$S\left\langle 2,3,5,6,8,9\right\rangle=
\left\langle 2,\underline{6},\underline{11},14\right\rangle$.

Now, for $\langle I;J\rangle\in M_3(n,q)$, where
$I=\left\langle i_1,i_2\dt i_r\right\rangle$, set
\[
S^{-1}\left\langle I;J\right\rangle=\left\langle x_1,x_2\dt x_r\mid y_1,y_2\dt y_r\right\rangle,
\]
where
\begin{gather*}
y_1=
\bcs
0 &\text{if \quad $i_1\not\in J$},\\
1 &\text{if \quad $i_1\in J$},
\ecs
\qquad
y_a=
\bcs
y_{a-1}+1 &\text{if \quad $i_a\not\in J$},\\
y_{a-1}+2 &\text{if \quad $i_a\in J$},
\ecs
\quad\text{for\quad$a\in[2,r]$},\\
x_a=i_a-y_a\quad\text{for\quad$a\in[1,r]$}.
\end{gather*}
Since $I$ is a 3--partition, $S^{-1}\left\langle I;J\right\rangle\in D(n,q)$
as it follows from Lemma \ref{lms}.

The definitions imply that the maps $S$ and $S^{-1}$ are mutually inverse.
Thus, the map $S$ is bijective.

\section{\bf Construction of a bijection $T:D(n,q)\rightarrow M_2(n,q)$}\label{sec03}

For any 1--partition $I=\langle i_1,i_2,\dots,i_q\rangle$, define
\[
\mu(I):=\left\{\,\min(r)\mid \langle i_r,i_{r+1}\dt i_q\rangle\;\text{is a $2$-partition}\,\right\}.
\]
We will construct partition $T(I)$ by induction on $\mu(I)$.
For $\mu(I)=1$, define $T(I)=I$.

\blr\label{lm1}
Let $I=\langle i_1,i_2,\dots,i_q\rangle\in D(n,q)$ and let $\mu(I)=k>1$.
Then there exists a unique
${s=s(I)\geqslant k}$ such that the marked partition
\[
\left\langle\,i_{k+1}-2,i_{k+2}-2\dt i_s-2,
\underline{i_{k-1}+i_k+2(s-k)},i_{s+1},i_{s+2},\dt i_q\,\right\rangle
\]
is a regular $2$-partition, where partitions
$\langle i_{k+1}-2,i_{k+2}-2\dt i_s-2\rangle$
and $\langle i_{s+1},i_{s+2},\dt i_q\rangle$are empty by definition
if $s=k$ and $s\geqslant q$, respectively.
\elr

\bp
Since $\langle i_k,i_{k+1},\dots,i_q\rangle$ is a $2$--partition,
it follows that
the sequence $i_s-2(s-k)$ does not decrease as $s\in\left[\,k,q\,\right]$ grows.
Therefore, there is a unique $s$ such that
\[
i_s-2(s-k)<i_{k-1}+i_k\leqslant i_{s+1}-2(s-k+1).
\]
This is equivalent to the required claim.
\ep

Let $T(I)$ be defined for all 1--partitions with $\mu(I)<k$, where $k\geqslant 1$.
Let $\mu(I)=k$ and $s=s(I)$. Set
\begin{gather*}
A(I):=
\big\langle\,i_1,i_2,\dots,i_{k-2},i_{k+1}-2,i_{k+2}-2,\dots,i_s-2\,\big\rangle,\\
B(I):=\big\langle\,\underline{i_{k-1}+i_k+2(s-k)},
i_{s+1},i_{s+2},\dots,i_q\,\big\rangle.
\end{gather*}
Since $\mu(A(I))<\mu(I)$ and $B(I)$ is a regular marked 2--partition,
then by inductive hypothesis the following marked partition is defined
\[
T(I):=T(A(I))\cup B(I).
\]
To complete the definition of $T$, it is sufficient to show that $T(I)\in M_2(n,q)$.

For brevity, set $l=\mu(A(I))$. For $l=1$, the claim follows from Lemma \ref{lm1}.
Let $l>1$. Then $s(A(I))>0$.
By definition, $s(A(I))\leqslant s-2$.

If this inequality is strict, then the required claim follows from Lemma \ref{lm1}.
Otherwise the marked part of partition $B(A(I))$
is equal to $i_{l-1}+i_l+2(s-l-2)$.
Since $l\leqslant k-2$, where $k=\mu(I)$, we have
\[
(i_{k-1}+i_k+2(s-k))-(i_{l-1}+i_l+2(s-l-2))
=(i_{k-1}-i_{l-1})+(i_k-i_l)-2(k-l-2)
\geqslant 2(k-l)\geqslant 4,
\]
because $i_{k-1}-i_{l-1}=i_k-i_l\geqslant 2(k-l-1)$.
Thus, $T(I)\in M_2(n,q)$.

Let us now construct the mapping $T^{-1}:M_2(n,q)\rightarrow D(n,q)$ inverse to $T$.

Define $T^{-1}\langle I;\emptyset\rangle=I$.
The next claim defines $T^{-1}\langle I;J\rangle$ for
$\langle I;J\rangle\in M_2(n,q)$ when $|J|=1$.

\blr\label{lm-t0}
For $\langle\,i_1,i_2,\dots,i_{a-1},\underline{i_a},i_{a+1},\dots,i_q\rangle\in M_2(n,q)$,
there is a unique $t=t(a)\in[1,a]$ such that
\begin{multline*}\label{exk}
T^{-1}\left\langle\,i_1,i_2,\dots,i_{a-1},\underline{i_a},i_{a+1},\dots,i_q\right\rangle:=\\[1mm]
\Big\langle\,i_1,i_2,\dots,i_{t-1},
\left\lfloor\frac{i_a}{2}\right\rfloor-(a-t),\,\left\lfloor\frac{i_a}{2}\right\rfloor-(a-t)+1,\,
i_t+2,i_{t+1}+2,\dots,i_{a-1}+2,i_{a+1},\dots,i_q\Big\rangle\in D(n,q),
\end{multline*}
where partitions $\langle i_1,i_2,\dots,i_{t-1}\rangle$ and
$\langle i_t+2,i_{t+1}+2,\dots,i_{a-1}+2\rangle$ are empty by definition if
$t=1$ and $t=a$, respectively.
\elr

\bp
If $a=1$ or $\left\lfloor\frac{i_a}{2}\right\rfloor>i_{a-1}$, set $t=a$.
Let $\left\lfloor\frac{i_a}{2}\right\rfloor\leqslant i_{a-1}$.
Since $\langle i_1,i_2,\dots,i_{a-1}\rangle$ is a 2--partition, then
\[
0\leqslant i_1-1<i_2-2<\dots<i_{a-1}-(a-1).
\]
Therefore, there is a minimal $t\in[1,a-1]$ such that
$\left\lfloor\frac{i_a}{2}\right\rfloor-a\leqslant i_t-t$.
This inequality is equivalent to the claim of Lemma.
\ep

For $\langle I;J\rangle\in M_2(n,q)$,
we will define $T^{-1}\langle I;J\rangle$ by induction on $|J|\geqslant 1$.
Assume that $T^{-1}\langle I;J\rangle$ is defined for all 2--partitions
with $|J|<k$ and, in addition, assume that
\beq\label{ineq}
\min\left(T^{-1}\langle I;J\rangle\right)\geqslant \left\lfloor\frac{i_a}{2}\right\rfloor-(a-1),
\qquad\text{where}\qquad i_a=\min(J).
\eeq
For $k=1$, this inequality is valid as it follows from Lemma \ref{lm-t0}.

Let $\langle I;J\rangle=\langle\,i_1,i_2,\dots,i_q\,;i_{a_1},i_{a_2},\dots,i_{a_k}\,\rangle\in M_2(n,q)$,
where $k\geqslant 2$, and let $t=t(a_1)$.
Set
\begin{gather*}
E\langle I;J\rangle:=
\left\langle\,i_1,i_2,\dots,i_{t-1},\,
\left\lfloor\frac{i_{a_1}}{2}\right\rfloor-(a_1-t),\,
\left\lfloor\frac{i_{a_1}}{2}\right\rfloor-(a_1-t)+1\right\rangle,\\
F\langle I;J\rangle:=
\left\langle\,i_t+2,i_{t+1}+2,\dots,i_{a_1-1}+2,i_{a_1+1},i_{a_1+2},\dots,i_q\;;
\;i_{a_2},i_{a_3},\dots,i_{a_k}\right\rangle\,.
\end{gather*}
The inductive hypothesis shows that the partition
\[
T^{-1}\langle I;J\rangle:=E\langle I;J\rangle\cup T^{-1}F\langle I;J\rangle
\]
is well defined.
To complete the induction step, it is sufficient to show that $T^{-1}\langle I;J\rangle\in D(n,q)$
and check inequality \eqref{ineq}.
But this inequality, obviously, follows from the required inclusion and Lemma \ref{lm-t0}.

By the inductive hypothesis we have
\[
\min\left(T^{-1}F\langle I;J\rangle\right)\geqslant \left\lfloor\frac{i_{a_2}}{2}\right\rfloor-(a_2-t-1).
\]
Therefore, to establish the inclusion $T^{-1}\langle I;J\rangle\in D(n,q)$
it suffices to show that
\beq\label{ineq1}
\left\lfloor\frac{i_{a_1}}{2}\right\rfloor-(a_1-t)+1<
\left\lfloor\frac{i_{a_2}}{2}\right\rfloor-(a_2-t-1),
\qquad\text{i.e.,}\qquad
\left\lfloor\frac{i_{a_2}}{2}\right\rfloor-\left\lfloor\frac{i_{a_1}}{2}\right\rfloor>a_2-a_1.
\eeq
The definition of regular marked 2--partition implies that
$i_{a_2}-i_{a_1}>2(a_2-a_1)$.
Since the numbers $i_1$ and $i_2$ are odd, the inequality \eqref{ineq1} follows.

A routine test shows that the maps $T$ and $T^{-1}$ are mutually inverse.
Thus, the map $T$ is bijective.
\medskip

For instance, from the definition of $T$ we obtain
\[
T:\;\langle 1,2,4,5,6,8\rangle\to\langle 1,2,4,8-2,\underline{5+6+2}\rangle=\langle 1,2,4,6,\underline{13}\rangle
\to\langle 4-2,6-2,\underline{1+2+2\cdot 2},\underline{13}\rangle.
\]
Therefore, $T\langle 1,2,4,5,6,8\rangle=\langle 2,4,\underline{7},\underline{13}\rangle$.

\brm\label{rem2}
At the price of making the arguments used to construct the mapping $T$
a bit more complicated one can construct a bijective mapping
$T_\lambda:D(n,q)\rightarrow M_\lambda(n,q)$ \textit{simultaneously}
for $\lambda=2$ and $\lambda=3$, where $T_2=T$.
We skip the precise definition of $T_3$ and just give an example of how it works:
\begin{multline*}
T_3:\;\langle 2,3,4,5,6,10,14\rangle\to\langle 2,3,4,11,10-3,14-3,\underline{5+6+2\cdot 3}\rangle=
\langle 2,3,4,7,11,\underline{17}\rangle\\
\to
\langle 2,7-3,11-3,\underline{3+4+2\cdot 3},\underline{17}\rangle=
\langle 2,4,8,\underline{13},\underline{17}\rangle\to
\langle 8-3,\underline{2+4+3},\underline{13},\underline{17}\rangle.
\end{multline*}
Therefore, $T_3\langle 2,3,4,5,6,10,14\rangle=\langle 5,\underline{9},\underline{13},\underline{17}\rangle$.
\erm

\section{\bf Marked $(\lambda,k)$-partitions}\label{sec04}

For a partition $I$, set $I^-=\min(I)$.

\bdr
A \emph{$(\lambda,k)$-partition} is a triple $(\lambda,k;I)$, where $k\geqslant 1$,
$I$ is a $\lambda$-partition, and $I^-\geqslant k$.
\edr

We write any $(\lambda,k)$-partition as a partition $I$, emphasising
that we treat $I$ as a $(\lambda,k)$-partition.
For example, $\langle 2,5,8\rangle$ considered as a $(3,1)$-partition, or as a $(3,2)$-partition,
or as a $(2,1)$-partition are different objects.

\bdr
A $(2,k)$-partition $I=\{i_1,i_2,\dots,i_q\}$ is called a \emph{special $(2,k)$-partition} whenever
\[
i_q<2(k+q-1).
\]
A $(3,k)$-partition $I=\{i_1,i_2,\dots,i_q\}$ is called a \emph{special $(3,k)$-partition} whenever
\begin{gather*}
i_q\leqslant 2k+3(q-1)\quad\text{if}\quad i_1>k,\\
i_q<2k+3(q-1)\quad\text{if}\quad i_1=k.
\end{gather*}
\edr

\brm
Let $S_{\lambda,k}(q)$ be the set of special $(\lambda,k)$-partitions of length $q$.
It is easy to show that
\[
|S_{2,k}(q)|=\binom{q+k-1}{k-1},\qquad\qquad
|S_{3,k}(q)|=\binom{q+k-1}{k-1}+\binom{q+k-2}{k-1}.
\]
\erm

\bdr
A $(\lambda,k)$-partition is called \emph{simple} if it is either special or dense.
\edr
\bdr
For any $(\lambda,k)$-partition $I$, there is a unique decomposition
${I=I_1\sqcup I_2\sqcup\dots\sqcup I_s}$,
where ${I_1,I_2,\dots,I_s}$ are simple  $(\lambda,k)$-partitions of the maximal possible length.
This decomposition is called the \emph{canonical form of $I$}.
The partitions $I_1,I_2,\dots,I_s$ are called the \emph{simple components of $I$}.
\edr

\bdr
Let $I$ be a $(\lambda,k)$-partition and let
$L\subset I$ be any non-special simple component of $I$
such that $\ind_\lambda(L)=1$. Then $L$ is called a
\emph{leading component of $I$}, and $L^-$
is called a \emph{leading part of $I$}.

The quantity of the leading components of
$I$ is called the \emph{index of $I$} and denoted by $\ind_{\lambda,k}(I)$.
\edr

For example, let $I=\langle 2,5,9,13,16\rangle$.
Then the canonical form of $I$ is
\[
I=
\bcs
\langle 2,5\rangle\sqcup\langle 9\rangle\sqcup\langle 13,16\rangle&\text{as a $(3,1)$-partition},\\
\langle 2,5,9\rangle\sqcup\langle 13,16\rangle&\text{as a $(3,2)$-partition}.
\ecs
\]
Therefore, $\ind_{3,1}(I)=2$ and $\ind_{3,2}(I)=1$.
Note also that $\ind_{2,1}(I)=3$ and $\ind_{2,1}(I)=2$.

\bdr
We say that a \emph{marked partition $\langle I;J\rangle$ is a
${(\lambda,k)}$-partition} if $I$ is a ${(\lambda,k)}$-partition;
we say that it is \emph{regular}
if $J$ is a subset of the set of leading parts of $I$.
\edr

In what follows we use the notation:
\[
\begin{array}{ll}
D_k(n,q)&\text{is the set of $1$-partitions of degree $n$, length $q$, and with minimal part $\geqslant k$}.\\[.5mm]
M_{\lambda,k}(n,q)&\text{is the set of the regular $(\lambda,k)$-partitions of degree $n$ and length $q$}.\\[.5mm]
\end{array}
\]

The next claim is the main result of this section.

\btr\label{t_maink}
For $\lambda=2$ or $\lambda=3$ and for any $k\geqslant 1$, we have
\beq\label{mainpk}
\prod_{r=k}^\infty(1+tx^r)=1+\sum_{q=1}^\infty\;
\sum_{n=\frac{q(q+1)}{2}+(k-1)q}^\infty\left|M_{\lambda,k}(n,q)\right|\,t^q\,x^n.
\eeq
In particular, $|D_k(n,q)|=|M_{\lambda,k}(n,q)|$.
\etr

\bp
For $k=1$, this is the result of Theorem \ref{t_main} since $M_{\lambda,1}(n,q)=M_\lambda(n,q)$.

Let $\langle I;J\rangle$ be a marked $(\lambda,k)$-partition, where
\beq\label{IJ}
I=\langle i_1,i_2,\dots,i_q\rangle,\qquad J=\langle i_{a_1},i_{a_2},\dots,i_{a_m}\rangle,
\eeq
and let ${I=I_1\sqcup I_2\sqcup\dots\sqcup I_s}$ be the canonical form of $I$.
Assume that for any $l\in[1,s]$, we have
\[
|J\cap I_l|=0\quad\text{if $I_l$ is not a leading component,
}\qquad\text{and}\qquad|J\cap I_l|\leqslant 1\quad\text{otherwise}.
\]
Then for any $r\in[1,m]$, there is a unique leading component $I_{t(r)}$ of $I$ such that
$i_{a_r}\in I_{t(r)}$.
Set
\[
\tau\langle I;J\rangle=\langle I;\tau(J)\rangle,\qquad\text{where}\qquad
\tau(J)=\{I^-_{t(1)},I^-_{t(2)},\dots,I^-_{t(m)}\}.
\]

Let us define a bijective map
$\delta:M_{\lambda,k}(n,q)\to M_{\lambda,k-1}(n-q,q)$ as follows.
For $\langle I;J\rangle\in M_{\lambda,k}(n,q)$, where for $I$ and $J$ notation \eqref{IJ} is used,
define the marked $(\lambda,k-1)$-partition $\sigma\langle I;J\rangle=\langle I^\prime;J^\prime\rangle$
by the formulas
\[
I^\prime=\left\langle i^\prime_1,i^\prime_2,\dots,i^\prime_q\right\rangle,
\quad\text{where}\quad
i_r^\prime=
\bcs
i_r-1&\text{if $i_r\not\in J$},\\
i_r-2&\text{if $i_r\in J$},
\ecs
\]
\[
J^\prime=\left\langle i^\prime_{a_1},i^\prime_{a_2},\dots,i^\prime_{a_m}\right\rangle,
\quad\text{where}\quad
i^\prime_{a_r}=i_{a_r}-2.
\]
Then the map $\tau$ is correctly defined on $\sigma\langle I;J\rangle$.
Now set $\delta\langle I;J\rangle:=\tau\sigma\langle I;J\rangle$.

A direct verification, which uses only the above definitions,
shows that ${\delta\langle I,J\rangle\in M_{\lambda,k-1}(n-q,q)}$.
It is easy to see that the map $\delta$ is invertible. Therefore, the map $\delta$ is bijective.
In particular,
\beq\label{M1q}
\left|M_{\lambda,k}(n,q)\right|=\left|M_{\lambda,k-1}(n-q,q)\right|=\dots=\left|M_{\lambda,1}(n-(k-1)q,q)\right|.
\eeq
Substituting $t\mapsto tx^{k-1}$ in formula \eqref{main1} and then applying formula \eqref{M1q}
we obtain formula \eqref{mainpk}.
\ep

\brm
Similarly to \S\ref{sec002}, formula $|D_k(n,q)|=|M_{\lambda,k}(n,q)|$ can be presented in the form
\[
\sum_{q=1}^\infty|D_k(n,q)|t^q=\sum_{q=1}^\infty\sum_{\alpha=0}^q|R_{\lambda,k}(n,q;\alpha)|\,(1+t)^\alpha t^q,
\]
where $R_{\lambda,k}(n,q;\alpha)$ denotes the set of $(\lambda,k)$-partitions
of degree $n$, of length $q$, and of index $\alpha$.

Using this formula and similar argumentation as in \S\ref{sec002} for $\lambda=2$, it is not difficult
to prove the following generalization of formula \eqref{Syl0}:
\beq\label{2k}
\prod_{r=\left\lceil\frac{k+1}{2}\right\rceil}^\infty(1-x^{2r-1})=
1+\sum_{q=1}^\infty(-1)^qx^{q(q+k-1)}\;\frac{{q+k-1\brack k-1}_x}{(1-x^{2k})(1-x^{2(k+1)})\dots(1-x^{2(q+k-1)})}\;,
\eeq
where
\[
{q+k-1\brack k-1}_x=\frac{(1-x^{q+1})(1-x^{q+2})\dots(1-x^{q+k-1})}
{(1-x)(1-x^2)\dots(1-x^{k-1})}
\]
is the Gaussian binomial coefficient (see \cite{MR1634067}, Ch.3).
\erm



\end{document}